%% file: paper.tex
\documentclass[10pt]{article}
\usepackage{amssymb,amsmath}
\usepackage{url}

\usepackage[hmargin = 1in, vmargin = 1in]{geometry}

\usepackage[parfill]{parskip}
\usepackage{graphicx}
\usepackage[font=bf]{caption}

\usepackage{subfig}
\usepackage{algorithm}
\usepackage{color}
\usepackage{algpseudocode}

\usepackage{fancyhdr}
\usepackage{bbm}
\usepackage{dsfont}

\newcommand\thefontsize[1]{{#1 The current font size is: \f@size pt\par}}

\newcommand{\R}{\mathbb{R}}
\newcommand{\N}{\mathbb{N}}
\newcommand{\probb}{\mathbb{P}}
\newcommand{\expv}{\mathbb{E}}
\renewcommand{\epsilon}{\varepsilon}
\renewcommand{\phi}{\varphi}
\newcommand{\prno}{\mathit{Pr}}
\newcommand{\reno}{\mathit{Re}}

\newcommand{\OM}{\Omega}
\newcommand{\TINT}{\int_{0}^T}

\newcommand{\sensor}[1]{#1^*} 

\newcommand{\abs}[1]{\left\lvert #1\right\rvert}

\newcommand{\norm}[1]{\left\lVert #1\right\rVert}

\newcommand{\set}[1]{\left\{ #1\right\}}

\newcommand{\paren}[1]{\left( #1\right)}

\let\p=\paren

\let\pb=\parenb

\newcommand{\sqparen}[1]{\left[ #1\right]}

\let\sp=\sqparen

\newcommand{\dprod}[2]{\left\langle #1, #2\right\rangle}
\newcommand{\dprodb}[2]{\bigl\langle #1, #2\bigr\rangle}
\newcommand{\dprodB}[2]{\Bigl\langle #1, #2\Bigr\rangle}

\newcommand{\pderiv}[3][]{\frac{\partial^{#1} #2}{\partial #3^{#1}}}

\newcommand{\diver}[2][]{\nabla_{\!#1} \cdot #2}

\newcommand{\grad}[2][]{\nabla_{\!#1} #2}
\newcommand{\lapl}[2][]{\triangle_{#1} #2}
\newcommand{\bdry}[1]{\partial #1}
\newcommand{\unit}[1]{\,{\rm #1}}

\newcommand{\uf}[1]{\textbf{DEPRECATED}}

\begin{document}

\begin{center}
\vspace{0.2in}
\noindent{\fontsize{14pt}{1em}\selectfont \textbf{
    Gradient-Based Estimation of Air Flow and Geometry Configurations in a Building Using Fluid Dynamic Adjoint Equations } }\\[14pt]

{\fontsize{11pt}{1.2em}\selectfont
  Runxin HE\textsuperscript{1*}, Humberto GONZALEZ\textsuperscript{2}
  \\[11pt]
  Department of Electrical \& Systems Engineering, Washington University in St.\ Louis,\\
  St.\ Louis, MO, USA \\[11pt]
  \textsuperscript{1}
  runxinhe@email.wustl.edu
  \\[11pt]
  \textsuperscript{2}
  hgonzale@wustl.edu
  \\[11pt]
  * Corresponding Author
  \\[33pt]
}
\end{center}

\input{sections/abstract}
\input{sections/intro}
\input{sections/problem_description}
\input{sections/gradient_algorithm}
\input{sections/exp_results}
\input{sections/conclusion}
\appendix
\input{sections/appendx}

\bibliographystyle{plain}
\bibliography{refs}

\end{document}

%% file: sections/abstract.tex
\section*{ABSTRACT}
Real-time estimations of temperature distributions and geometric configurations are important to energy efficient buildings and the development of smarter cities.
In this paper we formulate a gradient-based estimation algorithm capable of reconstructing the states of doors in a building, as well as its temperature distribution, based on a floor plan and a set of thermostats.
Our algorithm solves in real time a convection-diffusion Computer Fluid Dynamics (CFD) model for the air flow in the building as a function of its geometric configuration.
We formulate the estimation algorithm as an optimization problem, and we solve it by computing the adjoint equations of our CFD model, which we then use to obtain the gradients of the cost function with respect to the flow's temperature and door states.
We evaluate the performance of our method using simulations of a real apartment in the St.\ Louis area.
Our results show that the estimation method is both efficient and accurate, establishing its potential for the design of smarter control schemes in the operation of high-performance buildings.

%% file: sections/intro.tex
\section{INTRODUCTION}
\label{sec:intro}

Buildings currently account for more than 40\% of the total energy consumption in the U.S.~\cite{BEDB2011}, and they cost \$130 billion in energy leaks and inefficiencies~\cite{Council2014}.
For this reason many research groups have developed new control algorithms to improve the performance and efficiency of Heating, Ventilation, and Air Conditioning (HVAC) systems in buildings~\cite{wang2008supervisory, dounis2009advanced,afram2014theory,shaikh2014review}.

Among the many control algorithms used in smart building applications, Model Predictive Control~(MPC) stands out thanks to its flexible mathematical formulation, as well as accurate and robust responses in real-world implementations~\cite{huang2011model,xi2007support,morocsan2010building}.
Moreover, MPC has become the standard to solve complex constrained multivariate control problems in process control applications~\cite{afram2014theory}.
MPC has been used in HVAC control applications such as zoned temperature control~\cite{ma2012demand,he2015zoned} and overall temperature regularization~\cite{privara2011model} among others, experimentally showing significant increments in energy efficiency.

MPC algorithms require the use of dynamical models, which is used as predictors.
Compared to concentrated-parameter models such as the Resistance-Capacitance~(RC) networks~\cite{ma2012predictive}, Computational Fluid Dynamics~(CFD) models for building temperature control have a significant advantage since they naturally incorporate geometric and air flow information.
Moreover, CFD models can accurately describe short time scales, allowing us to reflect indoor climate changes in minute- or even second-level~\cite{jones1992computational,chow1996application}.
Hence, using MPC with CFD models as predictors enables the study of control and estimation strategies beyond the standard temperature control, such as controlling fan speeds or finding the optimal location for thermostats.

However, due to the predictive nature of its formulation, MPC is sensitive to mismatches in the dynamic prediction model and inaccurate initial state estimations, which might lead to steady-state offsets or even system instability~\cite{faanes2005offset}.
Therefore MPC algorithms are usually implemented in coordination with estimation algorithms capable of inferring relevant parameters and initial conditions from sensor data, such as Kalman filters.
Other optimization-based estimation algorithms have been developed in the past, such as the results in~\cite{isakov2006inverse} in Partial Differential Equation~(PDE) estimation.
Banks et al.\ viewed the parameters for the inverse problem as random variables, and used probabilistic inference methods to estimate the desired parameters~\cite{banks2001modelling, banks2011comparison}.
In~\cite{gutman1990identification} the authors fully discretized a weak form of Stokes Equations in time and space and identified the system's discontinuous parameters.

In this paper, we develop a gradient-based optimization method to estimate the doors configuration and temperature distribution in a building.
In particular, our contribution is twofold.
First, we mathematically formulate a gradient-based estimation method to identify real-time indoor climate distribution and a building's doors configuration based only on thermostatic data.
Second, we show the accuracy of our estimation method under a limited number of thermostats by simulating a real apartment in the St.\ Louis area.
Our results show that thermostatic information, when used together with CFD models, provide enough information to estimate most of the variables relevant for building climate control.
In other words, a handful of thermostats can provide information, such as the configuration of doors, without the need to physically install extra sensors in a building.

The paper is organized as follows.
The fluid dynamic model and finite element method are formulated in Section~\ref{sec:prb_dsc}.
We present the theoretical basis for our gradient-based estimation algorithm in Section~\ref{sec:grad_algo}.
Finally, our simulation results are presented in Section~\ref{sec:exm}. 

%% file: sections/problem_description.tex
\section{Computational Fluid Dynamic Model and Optimal Problem}
\label{sec:prb_dsc}

The kernel of our model is the incompressible Navier-Stokes equation, which is a good approximation for the coupling of temperature with free flow convection at atmospheric conditions~\cite{awbi1989,Dobrzynski2004}.
Throughout the paper we make two major simplifications to this model.
First, we assume that the air flow behaves as a laminar fluid which reaches steady-state behavior much faster than the temperature in the building.
Theoretical~\cite{awbi2003ventilation} and experimental~\cite{sun2003} results have shown that turbulent flows are present in residential building, such as in the area around HVAC vents, yet their overall effect in the temperature distribution is negligible.
Hence, we consider a stationary laminar Navier-Stokes equation to describe the fluid behavior, and a time-dependent equation to describe the temperature behavior.
Second, we consider only two-dimensional air flows moving parallel to the ground.
These assumptions reduce the accuracy of our model to some extent~\cite{van2013comparison}, yet they allow us to significantly simplify the computational complexity of our CFD-based control design.

Let $\Omega \subset \R^2$ be the area of interest, assumed to be bounded and connected, and let $\bdry{\Omega}$ be its boundary.
Let $u\colon \Omega \to \R^2$ be the \emph{stationary air flow velocity}, and $p\colon \Omega \to \R$ be the \emph{stationary air pressure} in $\Omega$.
Also, given $T > 0$, let $T_e\colon \Omega \times [0,T] \to \R$ be the \emph{temperature} in $\Omega$.
Then, following~\cite{landau_1971}, the non-dimensional temperature convection-diffusion model in $\Omega$ can be described by the following PDE:
\begin{equation}
  \label{eq:temp}
  \pderiv{T_e}{t}(x,t)
  - \grad[x] \cdot \pb{\kappa(x)\, \grad[x]{T_e}(x,t)}
  + u(x) \cdot \grad[x]{T_e}(x,t)
  = g_{T_e}(x,t),
\end{equation}
where $g_{T_e} \colon \Omega \times [0,T] \to \R$ represents the heat source in the room,
$\kappa:\Omega \to \R$ is the \emph{thermal diffusivity},
$\prno$ is the \emph{Prandtl number} of the air,
$\reno$ is the \emph{Reynolds number} of the air,
and $\grad[x]{} = \pb{\pderiv{}{x_1},\pderiv{}{x_2}}^T$ is the \emph{gradient operator}.
The initial condition of the temperature is:
\begin{equation}
  \label{eq:init_t}
  T_e(x,0) = \pi_0(x),\quad \text{for}\ x \in \Omega.
\end{equation}
Similarly, the non-dimensional stationary air flow in $\Omega$ is governed by the following incompressible Navier-Stokes stationary PDE:
\begin{gather}
  \label{eq:ns1}
  - \frac{1}{\reno} \lapl[x]{u}(x)
  + \pb{u(x) \cdot \grad[x]{}}\, u(x)
  + \grad[x]{p}(x)
  + \alpha(x)\, u(x)
  = g_u(x),\quad \text{and},\\
  \label{eq:ns2}
  \diver[x]{u}(x) = 0,
\end{gather}
where $g_u \colon \Omega \to \R^2$ represents all the external forces applied to the air (such as fans),
and $\lapl[x]{} = \pderiv[2]{}{x_1} + \pderiv[2]{}{x_2}$ is the \emph{Laplacian operator}.
The we introduce the \emph{viscous friction coefficient} $\alpha \colon \Omega \to \R$, following the technique in~\cite{gersborg2005topology}, to model different materials in $\Omega$.
Indeed, when the point $x$ corresponds to a material that blocks air, we choose $\alpha(x) \gg u(x)$, which results in $u(x) \approx 0$.
When the point $x$ corresponds to air, then we choose $\alpha(x) = 0$.

The building's exterior wall are denoted by $\Gamma_w \subset \bdry{\Omega}$, and the air inlet of the HVAC system is modeled as a gap in the wall, denoted $\Gamma_i \subset \bdry{\Omega}$.
Hence, $\Gamma_w \cup \Gamma_i = \bdry{\Omega}$.
The boundary condition for the temperature is:
\begin{equation}
  \label{eq:bnd_t}
  T_e(x) \equiv T_A,\quad \text{for}\ x \in \bdry{\Omega},
\end{equation}
where $T_A$ is the \emph{atmospheric temperature}.
We only apply a boundary condition for the pressure equation at the inlet, setting $p(x) \equiv p_A$ for each $x \in \Gamma_i$, where $p_A$ is the \emph{atmospheric pressure}.
We do not define a boundary condition for the air flow at the inlet $\Gamma_i$, and we set the air flow at the exterior wall as follows:
\begin{equation}
  \label{eq:bnd_u1}
  u(x) \equiv 0,\quad \text{for}\ x \in \Gamma_w.
\end{equation}

We assume that there are $n_t$ thermostats in the building.
The $i$-th thermostat is located at $x_i \in \Omega$, and samples the temperature in a neighborhood averaged using the \emph{bump weight function} $\Phi_i(x) = \sigma\, \exp\pb{-\p{r^2 - \|x - x_i\|^2}^{-1}}$ for $\norm{x - x_i} < r$, and $\Phi_i(x) = 0$ otherwise, where $\sigma > 0$ is a normalization factor such that $\int_{\Omega} \Phi_i(x)\, \mathrm{d}x = 1$.

We also assume that there are $n_d$ doors in the building.
We define $\theta_i \in \set{0,1}$ as the configuration of the $i$-th door, i.e., $\theta_i = 1$ when the $i$-th door is open, and $\theta_i = 0$ when is closed.
Let $\Omega_{\theta_i} \subset \Omega$ be the area occupied by the $i$-th door when it is closed, and let $\mathbb{I}_i$ be the indicator function of $\Omega_{\theta_i}$, i.e., $\mathbb{I}_i(x) = 1$ for $x \in \Omega_{\theta_i}$, and $\mathbb{I}_i(x) = 0$ otherwise.

When the door configuration changes, so does the prediction generated by our CFD model in equations~\eqref{eq:temp} and~\eqref{eq:ns1}.
In particular, the parameters $\alpha$ and $\kappa$ change for each $x \in \Omega_{\theta_i}$ as a function of $\theta_i$.
We model this relation by defining $\alpha\colon \Omega \times \set{0,1}^{n_d} \to \R$ and $\kappa\colon \Omega \times \set{0,1}^{n_d} \to \R$ as follows:
\begin{equation}
  \label{eq:kappa}
  \alpha(x,\theta)
  = \alpha_0 + \sum_{i=0}^{n_d} \p{1-\theta_i}\, \p{\alpha_w-\alpha_0}\, \mathbb{I}_i(x),
  \quad \text{and} \quad
  \kappa(x,\theta)
  = \kappa_0 + \sum_{i=0}^{n_d} \p{1-\theta_i}\, \p{\kappa_w-\kappa_0}\, \mathbb{I}_i(x),
\end{equation}
where $\alpha_0$ and $\kappa_0$ are the parameters for open air, while $\alpha_w$ and $\kappa_w$ are the parameters for solid walls.
Note that both $\alpha$ and $\kappa$ are affine functions of $\theta \in \R^{n_d}$.

Now, using binary values for each $\theta_i$ means that our estimation algorithm will have to use combinatorial methods, which tend to scale poorly in both computation time and computational resources.
To avoid this problem we relax the binary parameters $\theta_i \in \set{0,1}$, instead allowing them to belong to the unit interval $[0,1]$.
Although for each $\theta_i$ only the extreme values have meaningful physical interpretations, non-integer values can theoretically be interpreted as averaged observations over the optimization horizon, as explained in~\cite{Vasudevan2013a,Vasudevan2013b}.
For example, if throughout the optimization horizon a door is open half the time, and closed half the time, it is likely that we will observe $\theta_i \approx 0.5$.
The relaxation of each $\theta_i$ is also important in our numerical calculations, since it transforms the optimization program from a mixed-integer program to a more convenient nonlinear format~\cite{gersborg2005topology}.

Now we can formulate our main estimation algorithm to compute the door configuration $\theta$ and the initial temperature $\pi_0$ using the information from the $n_t$ thermostats in the building.
Given an arbitrary estimation time horizon, say $[0,T]$, we write our optimal estimation problem as follows:
\begin{equation}
  \begin{aligned}
    \label{eq:cost}
    \min_{\pi_0\colon \OM \to \R,\, \theta \in \R^{n_d}}
    &J\p{\pi_0,\theta} =
    \sum_{i=1}^{n_t} \TINT \p{\int_{\Omega_i} \Phi_i\, T_e(x,t;\pi_0,\theta)\, \mathrm{d}x - \sensor{T}_{e,i}}^2 \mathrm{d}t
    + \eta_0\, \sum_{i=1}^{n_t} \p{\int_{\Omega_i} \Phi_i\, \pi_0\, \mathrm{d}x - \sensor{\pi}_{0,i}}^2
    + \eta_1\, \norm{\pi_0}^2_{\OM}, \\
    \text{subject to:}\ \
    &\text{partial differential equations~\eqref{eq:temp}, \eqref{eq:ns1}, and~\eqref{eq:ns2},} \\
    &\text{boundary and initial conditions~\eqref{eq:init_t}, \eqref{eq:bnd_t}, and~\eqref{eq:bnd_u1},} \\
    &0 \leq \theta_i \leq 1, \quad \forall i \in \set{1,\dotsc,n_d},
  \end{aligned}
\end{equation}
were, $\eta_0, \eta_1 > 0$ are weight parameters, $T_e(x,t;\pi_0,\theta)$ is the unique solution of equation~\eqref{eq:temp} with initial condition $\pi_0$ and configuration $\theta$, $\sensor{T}_{e,i}(t)$ is the time signal obtained from the $i$-th thermostat over the horizon $[0,T]$, and $\sensor{\pi}_{0,i}$ is just notation for the initial thermostat temperature, i.e., $\sensor{\pi}_{0,i} = \sensor{T}_{e,i}(0)$.

%% file: sections/gradient_algorithm.tex
\section{ADJOINT-BASED GRADIENT COMPUTATION}
\label{sec:grad_algo}

In this section we develop a numerical algorithm to solve the optimization problem defined in equation~\eqref{eq:cost}.
We use a gradient-based optimization algorithm to find local minimizers of our optimization problem, where the gradients are computed using the adjoint equations of the CFD model, similar to the techniques in~\cite{gunzburger2000adjoint} and~\cite{yang2014reaction}.
We then discretize the adjoint equations using the Finite Element Method~(FEM), resulting in a practical algorithm which we test in Section~\ref{sec:exm}.

\subsection{Adjoint Equations and Fr\'echet Derivatives}
\label{sec:adj_eq}

In order to derive our CFD model's adjoint equations, first we need to write the \emph{Lagrangian function} of the optimization problem~\cite{giles1997adjoint, gunzburger2000adjoint}.
Let $\set{\lambda_i}_{i=1}^6$ be the set of \emph{Lagrange multipliers}, or \emph{adjoint variables}, each associated to one of the equations~\eqref{eq:temp} to~\eqref{eq:bnd_u1} and defined in its respective dual space.
Then, the Lagrangian function of our optimal estimation problem is:
\begin{multline}
  \label{eq:lg}
    L\p{T_e,u,p,\pi_0,\theta,\set{\lambda_i}_{i=1}^6}
    = J\p{\pi_0,\theta}
    + \dprodB{\lambda_1}{
      \pderiv{T_e}{t}
      - \grad[x] \cdot \p{\kappa(x)\, \grad[x]{T_e}}
      + u \cdot \grad[x]{T_e}
      - g_{T_e}}_{\OM \times [0,T]}
    + \dprod{\lambda_4}{T_e}_{\partial \OM \times [0,T]}
    + \\
    + \dprodB{\lambda_2}{
      - \frac{1}{\reno} \lapl[x]{u}
      + \p{u \cdot \grad[x]{}} u
      + \grad[x]{p}
      + \alpha\, u
      - g_u}_{\OM}
    + \dprod{\lambda_3}{\diver[x]{u}}_{\OM}
    + \dprod{\lambda_5}{u}_{\Gamma_w}
    + \dprod{\lambda_6}{T_e(0,\cdot) - \pi_0}_{\OM},
\end{multline}
where $\dprod{f_1}{f_2}_S = \int_S f_1(z)\, f_2(z)\, \mathrm{d}z$ is the inner product of the Hilbert space of square integrable functions ${\cal L}^2(S)$.
We write the necessary conditions for optimality using Galerkin methods~\cite{girault1979finite}, i.e., by setting the inner product of the partial derivatives of $L$ with respect to all the dual directions equal to zero.
That is, we look for solutions such that $\dprodb{\pderiv{L}{T_e}}{w}_{\OM \times [0,T]} = 0$, $\dprodb{\pderiv{L}{u}}{v}_{\OM} = 0$, and $\dprodb{\pderiv{L}{p}}{q}_{\OM} = 0$ for each set of functions $(w,v,q)$ in the respective dual spaces, and sufficiently weakly differentiable.
As detailed in Appendix~\ref{sec:appendix_a}, the conditions above are satisfied when the dual variables satisfy:
\begin{gather}
  \label{eq:ad_pde1}
  - 2\, \sum_{i=1}^{n_t} \p{\int_{\Omega_i} \Phi_i(z)\, T_e(z,t)\, \mathrm{d}z - \sensor{T}_{e,i}(t)}
  + \pderiv{\lambda_1}{t}(x,t)
  + \grad[x]{} \cdot \p{\kappa(x)\, \grad[x]{\lambda_1}(x,t)}
  + u(x) \cdot \grad[x]{\lambda_1}(x,t)
  = 0, \\
  \lambda_6(x) = \lambda_1\p{x,0},\\
  \label{eq:ad_pde2}
  \TINT \lambda_1(x,t)\, \grad[x]{T_e}(x,t)\, \mathrm{d}t
  + \alpha(x)\, \lambda_2(x)
  - \frac{1}{\reno} \lapl[x]{\lambda_2}(x)
  - u(x) \cdot \grad[x]{\lambda_2}(x)
  + \lambda_2(x) \cdot \grad[x]{u}(x)
  - \grad[x]{\lambda_3}(x)
  = 0,\quad \text{and},\\
  \label{eq:ad_pde3}
  \diver[x]{\lambda_2}(x)
  = 0,
\end{gather}
with boundary conditions $\lambda_1\p{x,t} = 0$ and $\lambda_2\p{x,t} = 0$ for each $x \in \bdry{\Omega}$ and $t \in [0,T]$, together with final condition $\lambda_1\p{x,T} = 0$ for each $x \in \Omega$.
The adjoint functions $\lambda_4$ and $\lambda_5$ are irrelevant to our Fr\'echet derivative calculation, therefore we omit them from this presentation.

Now we can compute the Fr\'echet derivatives of the cost function with respect to $\theta$ and $\pi_0$.
Consider a parameter change from $\p{\theta, \pi_0}$ to $\p{\theta+\delta \theta, \pi_0 + \delta \pi_0}$.
Since both $\alpha$ and $\kappa$ are affine in $\theta$, these variations will result in changes from $\p{\alpha, \kappa}$ to $\p{\alpha + \delta \alpha, \kappa + \delta \kappa}$, which will also imply changes from $\p{T_e,u,p}$ to $\p{T_e+\delta T_e, u + \delta u, p+\delta p}$.
As detailed in Appendix~\ref{sec:appendix_b}, these variations allow us to compute a first-order approximation of the cost function $J$, which result in:
\begin{equation}
  \label{eq:alpha_kappa_deriv}
  \dprodb{\mathcal{D}_{\alpha}J}{\delta \alpha}_{\OM}
  = \dprod{\lambda_2 \cdot u}{\delta \alpha}_{\OM},\quad \text{and} \quad
  \dprodb{\mathcal{D}_{\kappa}J}{\delta \kappa}_{\OM}
  = \TINT \dprodb{\grad[x]{\lambda_1} \cdot \grad[x]{T_e}}{\delta \kappa}_{\OM}\, \mathrm{d}t,
\end{equation}
and using the chain rule and the formulas in equation~\eqref{eq:alpha_kappa_deriv} we get the desired directional derivatives for $J$:
\begin{equation}
  \label{eq:theta_i0_deriv}
  \dprodb{\mathcal{D}_{\pi_0}J}{\delta \pi_0}_{\OM}
  = \dprod{\grad[\pi_0]{J} - \lambda_6}{\delta \pi_0}_{\OM},
  \quad \text{and} \quad
  \mathcal{D}_{\theta}J \cdot \delta \theta
  = \sum_{i=1}^{n_d} \p{\dprodB{\mathcal{D}_{\alpha}J}{\pderiv{\alpha}{\theta_i}}_{\OM} + \dprodB{\mathcal{D}_{\kappa}J}{\pderiv{\kappa}{\theta}}}\, \delta\theta_i.
\end{equation}
Note that both directional derivatives are linear bounded operators, hence they are also Fr\'echet derivatives as desired.

\subsection{Gradient-Based Optimization Algorithm}
\label{sec:grad_method}

Using the closed-form formulas for the Fr\'echet derivatives of $J$ with respect to $\pi_0$ and $\theta$, we build a gradient-based optimization algorithm to solve the problem in equation~\eqref{eq:cost} using a projected-gradient method~\cite[Chapter~18.6]{nocedal2006numerical}.

First, we find descent directions $\delta\pi_0$ and $\delta\theta$ as solutions of the following Quadratic Program~(QP) with value $V$:
\begin{equation}
  \label{eq:qp_grad}
  \begin{aligned}
    V =
    \min_{\delta\pi_0\colon \Omega \to \R,\, \delta\theta \in \R^{n_d}}
    &\dprod{\mathcal{D}_{\pi_0}J}{\delta\pi_0}_{\Omega}
    + \mathcal{D}_{\theta}J \cdot \delta\theta
    + \frac{\gamma}{2}\, \norm{\delta\pi_0}_\Omega^2
    + \frac{\gamma}{2}\, \norm{\delta\theta}^2,\\
    \text{subject to:}\ \
    &0 \leq \theta_i + \delta\theta_i \leq 1, \ \forall i \in \set{1,\dotsc,n_d},
  \end{aligned}
\end{equation}
where $\gamma > 0$ is a parameter.
The QP in equation~\eqref{eq:qp_grad} is derived using first-order approximations for the cost function using the derivatives in equation~\eqref{eq:theta_i0_deriv}, together with a condition to guarantee the feasibility of the desired direction.
Note that $V \leq 0$, since $\delta\pi_0 = 0$ and $\delta\theta = 0$ always belong to the feasible set.
Hence, if $V = 0$ then our method cannot find further descent directions, and it thus terminates.

Second, a step size is computed using the following Armijo line search method:
\begin{equation}
  \label{eq:armijo}
  \begin{aligned}
  \beta = &\arg\max_{j \in \N} \bar{\beta}^j,\\
  &\text{subject to:}\ \
  J\pb{\pi_0 + \bar\beta^j\, \delta\pi_0, \theta + \bar\beta^j\, \delta\theta} - J(T_e,\pi_0,\theta) \leq \bar\alpha\, \bar\beta^j\, V.
  \end{aligned}
\end{equation}
where $\bar\alpha, \bar\beta \in (0,1)$ are parameters.

\begin{algorithm}[tp]
  \caption{Gradient-based estimation algorithm}
  \label{alg:grd_est}
  \begin{algorithmic}[1]
    \Require Initial values for $\theta$ and $\pi_0$.
    \Loop
    \State \label{step:ns}
    Compute $T_e$, $u$, and $p$ by solving the CFD model in equations~\eqref{eq:temp} to~\eqref{eq:bnd_u1}.
    \State \label{step:adj}
    Compute $\lambda_1$, $\lambda_2$, $\lambda_3$, and $\lambda_6$ by solving the adjoint equations~\eqref{eq:ad_pde1} to~\eqref{eq:ad_pde3}.
    \State \label{step:grad}
    Compute the gradients $\mathcal{D}_{\pi_0}$ and $\mathcal{D}_{\theta}$ in equation~\eqref{eq:theta_i0_deriv}.
    \State Compute the projected-gradient descent directions $(\delta\pi_0, \delta\theta)$ by solving the QP in equation~\eqref{eq:qp_grad}, with value $V$.
    \If{$V = 0$}
    \State Stop.
    \EndIf
    \State Compute the step size $\beta$ using the Armijo line search method in equation~\eqref{eq:armijo}.
    \State Update $\pi_0 \leftarrow \pi_0 + \delta\pi_0$ and $\theta \leftarrow \theta + \beta\, \delta\theta$.
    \EndLoop
  \end{algorithmic}
\end{algorithm}

Our gradient-based optimization method is detailed in Algorithm~\ref{alg:grd_est}.
Steps~\ref{step:ns} and~\ref{step:adj}, are numerically solved using FEM discretizations, implemented using the FEniCS package~\cite{logg2012automated}.

%% file: sections/exp_results.tex
\section{EXPERIMENTAL RESULTS}
\label{sec:exm}

\begin{figure}[tp]
  \centering
  \begin{minipage}[c]{.4\linewidth}
    \includegraphics[width=\linewidth]{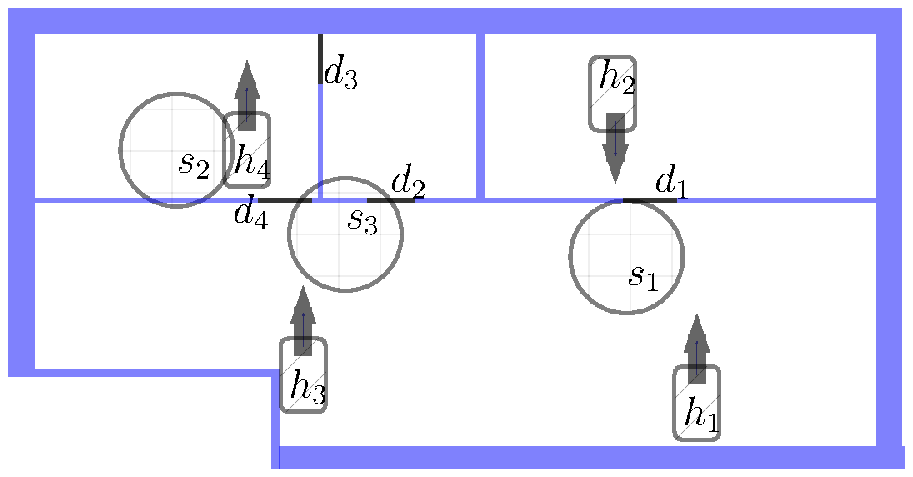}
    \caption{Floor plan of the apartment simulated in Sec.~\ref{sec:exm}.}
    \label{fig:build}
  \end{minipage}%
  \hspace{2em}%
  \begin{minipage}[c]{.55\linewidth}
  \subfloat[%
    Average estimation error of $\theta$.
  ]{%
    \label{fig:exp1_dr}%
    \hspace{1.5em}%
    \includegraphics[height=1.8in,trim=0 10 10 20,clip]{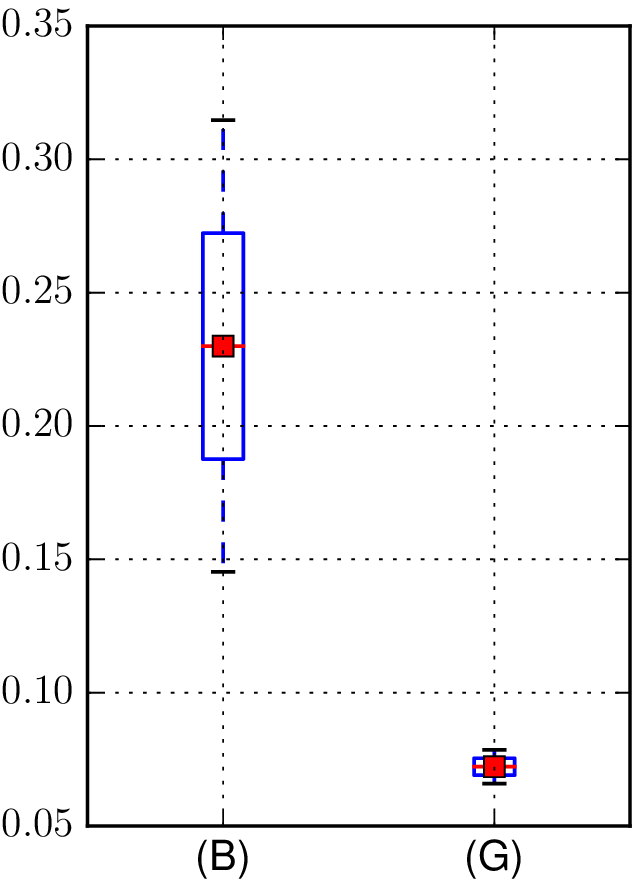}%
    \hspace{1.3em}%
  }%
  \hfill%
  \subfloat[%
    Relative estimation error of $\pi_0$.
  ]{%
    \label{fig:exp1_tem}%
    \hspace{1.5em}%
    \includegraphics[height=1.8in,trim=0 10 10 20,clip]{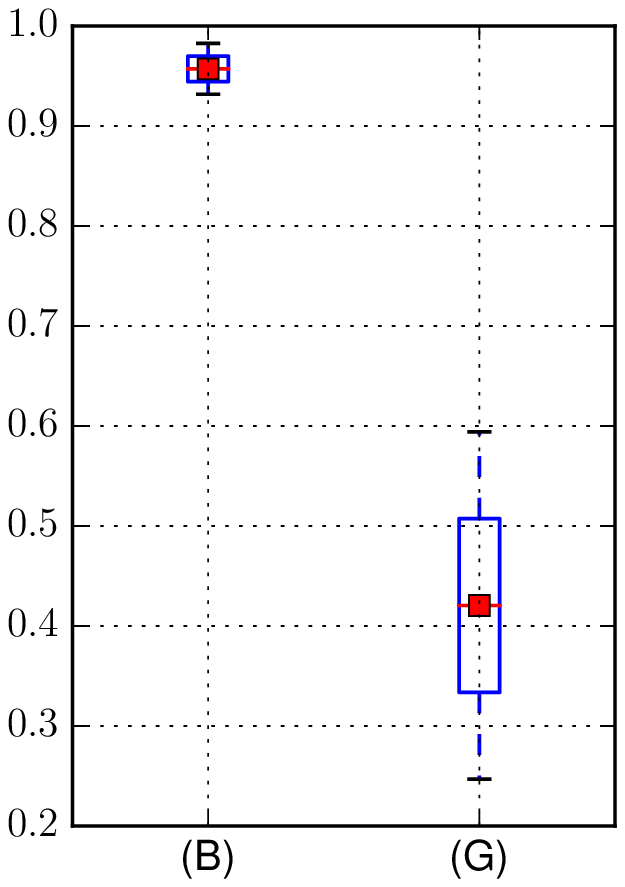}%
    \hspace{1.3em}%
  }%
  \caption{%
    Results of the experiments in Sec.~\ref{subsec:exp1}.
    Columns:
    (B)~Banks and Bihari's method,
    (G)~Algorithm~\ref{alg:grd_est}.
  }
  \label{fig:exp1_stats}
  \end{minipage}
\end{figure}

We applied our estimation algorithm to a simulated St.\ Louis area apartment with $n_d=4$ doors, labeled $\set{d_i}_{i=1}^4$, and $n_t=3$ thermostats, labeled $\set{s_i}_{i=1}^3$.
The floor plan of the apartment is shown in Figure~\ref{fig:build}, with dimensions $7.6 \times 16.8 \unit{m^2}$ (approx.\ $1375 \unit{sq\ ft}$).
The apartment is equipped with 4~HVAC vents, labeled $\set{h_i}_{i=1}^4$.
We assume that each vent is endowed with a fan acting on a $1 \times 0.5 \unit{m^2}$ area, and oriented in a fix direction.

The CFD model is governed by the constants $\reno = 10^2$, $\alpha_0 = 0$ and $\kappa_0 = 10^{-2}$ when $x \in \OM$ corresponds to free air, while $\alpha_w = 10^3$ and $\kappa_w = 10^{-4}$ when $x \in \OM$ corresponds to a wall.
The atmospheric pressure is $p_A = 101.3 \unit{kPa}$, and the atmospheric temperature is $T_A = 23.83 \unit{^\circ C}$.
We assume that $h_1$ and $h_2$ work at a low output setting, producing $0.1 \unit{kW}$ of heat and an air flow speed of $0.1 \unit{m/s}$.
On the other hand, $h_3$ and $h_4$ work at a normal setting, producing $4 \unit{kW}$ of heat and and an air flow speed of $0.5 \unit{m/s}$.
The time horizon is $300 \unit{s}$, sampled uniformly at $10 \unit{s}$ steps.
The sensors' observation radius is $r = 1.0 \unit{m}$.
The parameters in~\eqref{eq:cost} are set to $\eta_0 = 1.0$, $\eta_1 = 0.1$.
The parameter in~\eqref{eq:qp_grad} is set to $\gamma = 1.0$.
The parameters in~\eqref{eq:armijo} are set to $\bar{\alpha} = 0.01$ and $\bar{\beta} = 0.7$.
We wrote our code in Python, the FEM discretization was computed using tools from the FEniCS Project~\cite{logg2012automated}, and the building plan was discretized into $n_{\text{elem}} = 6276$ elements.

\subsection{Probabilistic Estimation Method}
\label{subsec:comp_method}

In our experiments below we compare our estimation method with a probabilistic-based estimation algorithms formulated in~\cite{banks2001modelling}, and applied to problems involving parameter estimation of differential equations~\cite{banks2011comparison}.
Under Banks and Bihari's framework, $\pi_0$ and $\theta$ are random variables with unknown probability distributions, thus the estimation problem is formulated such that we aim to find the optimal distributions that would most likely produce the acquired sensor data in expectation.
Due to space constraints we omit a detailed description of this method, we refer interest readers to~\cite{banks2001modelling}.

Let $\pi_{0,\Delta}$ be the FEM discretization of $\pi_0$, hence $\pi_{0,\Delta} \in \R^{n_{\text{elem}}}$.
We assume that $\theta$ and $\pi_{0,\Delta}$ follow probability distributions $\probb(\theta)$ and $\probb(\pi_0)$.
In the particular case of $\theta$, since it is a vector of independent binary variables, its distribution is $\probb(\theta) = \prod_{i=1}^{n_d} p_i^{\theta_i}\, (1-p_i)^{1 - \theta_i}$.
We assume that $\theta$ and $\pi_{0,\Delta}$ are independent.

Banks and Bihari's estimation algorithm relies on closed-form formulas of the expected values of each of the random variables in the cost function.
Using the cost function in equation~\eqref{eq:cost}, the only nontrivial expected value is that of $T_{e,\Delta}(x,t;\pi_{0,\Delta},\theta)$, the FEM discretization of $T_e$.
Note that given $x \in \Omega$, $t \in [0,T]$, and $\theta \in \set{0,1}^{n_d}$, then $T_{e,\Delta}(x,t;\pi_{0,\Delta},\theta)$ is a linear function of $\pi_{0,\Delta}$, hence, as shown in~\cite[Chapter~3]{Kumar2015}, the conditional expected value of $T_{e,\Delta}$ is $\expv\sp{T_e(x,t;\pi_{0,\Delta},\theta) \mid \theta} = T_{e,\Delta}(x,t;\expv\sp{\pi_{0,\Delta}},\theta)$ for each pair $(x,t)$.
Then, using Bayes' rule:
\begin{equation}
  \label{eq:exp_temp_2}
  \expv\sp{T_{e,\Delta}(x,t;\pi_{0,\Delta},\theta)} = \sum_{\theta \in \set{0,1}^{n_d}} T_{e,\Delta}(x,t;\expv\sp{\pi_{0,\Delta}},\theta)\, \probb(\theta),\quad \forall x \in \Omega,\ t \in [0,T].
\end{equation}
It is worth noting that the cardinality of $\set{0,1}^{n_d}$ is $2^{n_d}$, hence each evaluation of equation~\eqref{eq:exp_temp_2} involves solving a PDE an exponentially growing number of times as a function of $n_d$.

\subsection{Estimation Using Three Thermostats}
\label{subsec:exp1}

We run both estimation algorithms, Banks and Bihari's method and Algorithm~\ref{alg:grd_est}, under 6~different combinations for $\theta \in \set{(0,0,0,0), (0,0,1,1), (0,1,1,0), (0,1,1,1), (1,0,0,0), (1,1,1,1)}$ and 2~different initial temperatures $\pi_0$.
Since Algorithm~\ref{alg:grd_est} converges to local minimizers, we also run 5~estimations for each pair $(\pi_0,\theta)$ initializing the algorithm with different values.

\begin{figure}[tp]
  \centering
  \subfloat[%
  Initial temperature $\pi_0$ in $\unit{^\circ C}$ with respect to $T_A$.
  ]{%
    \label{fig:exp1_temp_iniB}%
    \includegraphics[width=.32\linewidth,trim=30 0 15 0,clip]{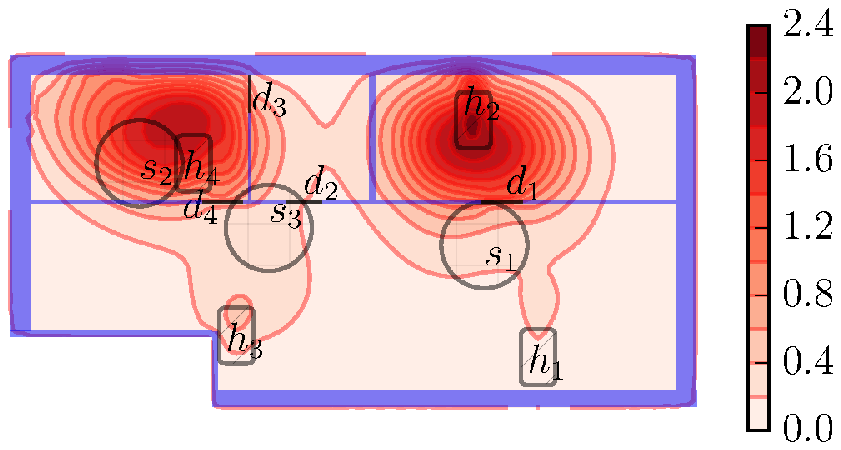}%
  }%
  \hfill%
  \subfloat[%
  Estimation error of $\pi_0$ by Banks and Bihari's method.%
  ]{%
    \label{fig:exp1_tem_bB}%
    \includegraphics[width=.32\linewidth,trim=30 0 15 0,clip]{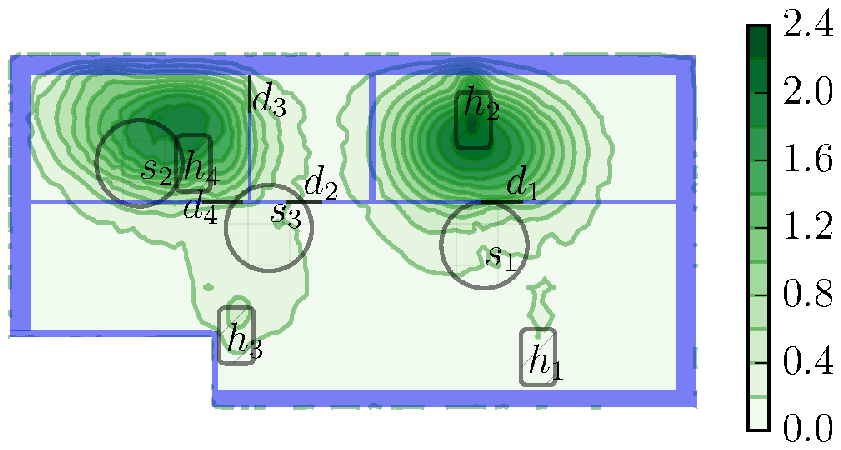}%
  }%
  \hfill%
  \subfloat[%
  Estimation error of $\pi_0$ by Algorithm~\ref{alg:grd_est}.%
  ]{%
    \label{fig:exp1_temp_iB}%
    \includegraphics[width=.32\linewidth,trim=30 0 15 0,clip]{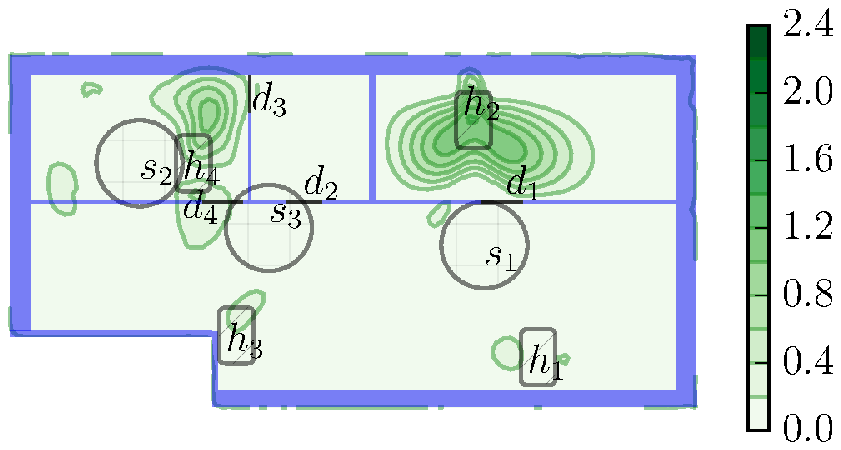}%
  }%
  \caption{%
    Results of the experiments in Sec.~\ref{subsec:exp1}.
    All doors are closed.
  }
  \label{fig:exp1_temp}
\end{figure}

Figure~\ref{fig:exp1_dr} shows a bar plot of the average estimation errors of $\theta$, calculated as $e_\theta = \frac{1}{n_d} \sum_{i=1}^{n_d} \abs{\theta_i - \hat{\theta}_i}$, where $\hat{\theta}$ is either the estimated probability distribution from Banks and Bihari's method, or the estimated relaxed configuration from Algorithm~\ref{alg:grd_est}.
Figure~\ref{fig:exp1_tem} shows a similar bar plot for the relative estimation error of $\pi_0$, calculated as $e_{\pi_0} = \frac{\norm{\pi_0 - \hat{\pi}_0}_\Omega}{\norm{\pi_0}_\Omega}$, where $\hat{\pi}_0$ is either the estimated expected value of $\pi_0$ from Banks and Bihari's method, or the estimated initial distribution from Algorithm~\ref{alg:grd_est}.
From these results we can observe that Algorithm~\ref{alg:grd_est} is significantly more accurate than the probabilistic method in estimating both door configuration and initial temperature.
It is worth noting that the average error of Algorithm~\ref{alg:grd_est} in Figure~\ref{fig:exp1_dr} is small enough so that one can use a constant threshold to convert from relaxed values of $\theta$ to binary values.
Also, the accuracy of our results enables further smart applications, such as the locating the residents in a building just by using thermostat data.
In Figure~\ref{fig:exp1_temp} we show the actual initial temperature distribution for one configuration~$\theta$, and the estimation errors by both algorithms.
These results show that even with the temperature of three points, Algorithm~\ref{alg:grd_est} can accurately reconstruct the initial temperature distribution in the building, thus enabling advanced control methods such as MPC to significantly improve the energy efficiency of the HVAC system~\cite{he2015zoned}.

\subsection{Estimation Using One Thermostat}
\label{subsec:exp2}

Now we only assume that only one thermostat, $s_1$, is functional.
The motivation for this experiment is to show the performance of both estimation algorithms in a realistic scenario, since most residential buildings' HVAC systems operate using a single thermostat.
We simulated the same scenarios as in Section~\ref{subsec:exp1}.

\begin{figure}[tp]
  \centering
  \begin{minipage}[c]{.65\linewidth}
  \subfloat[%
  Average estimation error of $\theta$.
  ]{%
    \label{fig:exp2_dr}%
    \includegraphics[height=1.8in,trim=0 10 10 20,clip]{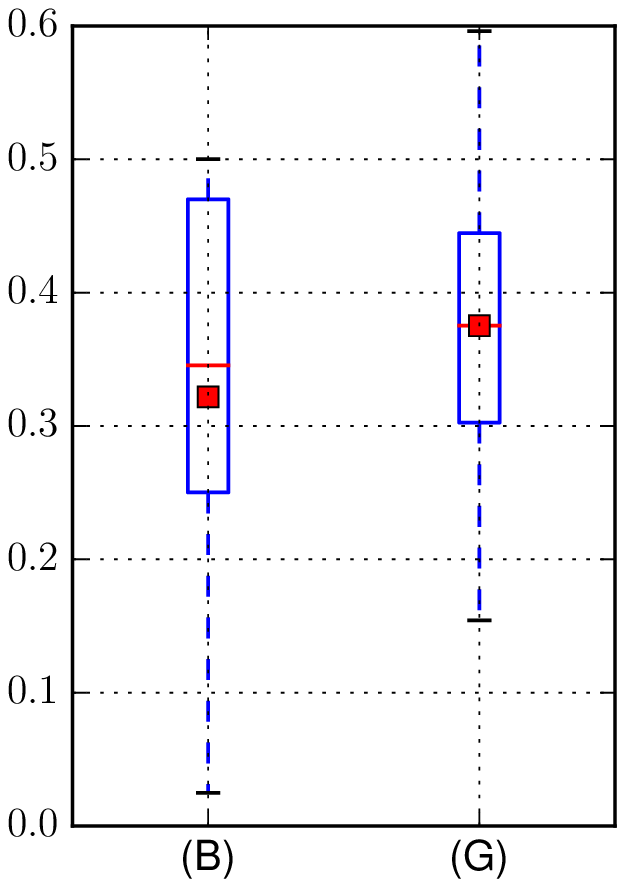}%
  }%
  \hfill%
  \subfloat[%
  Estimation error of $\theta_1$.
  ]{%
    \label{fig:exp2_dr1}%
    \includegraphics[height=1.8in,trim=0 10 10 20,clip]{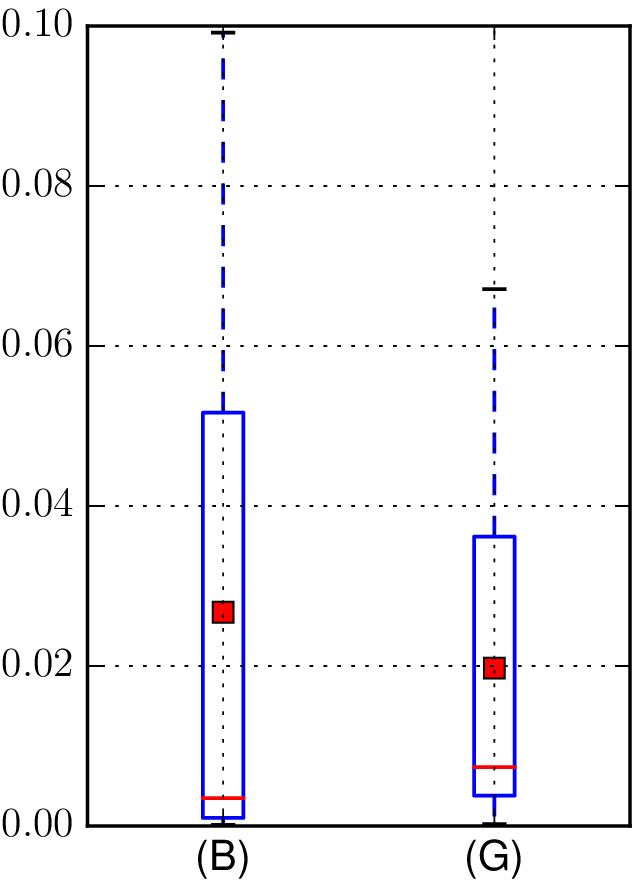}%
  }%
  \hfill%
  \subfloat[%
  Relative estimation error of $\pi_0$.
  ]{%
    \label{fig:exp2_tem}%
    \includegraphics[height=1.8in,trim=0 10 10 20,clip]{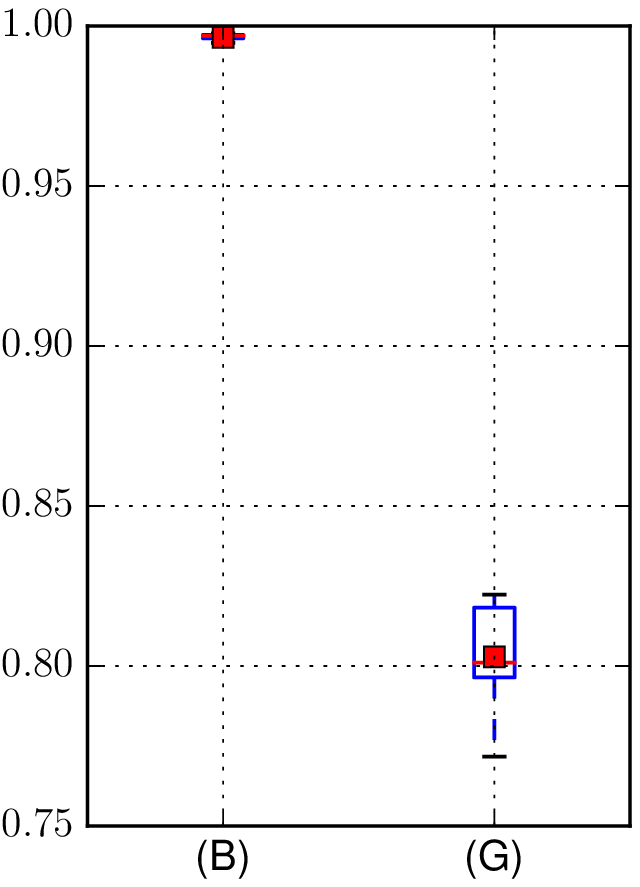}%
  }%
  \caption{%
    Results of the experiments in Sec.~\ref{subsec:exp2}.
    Columns:
    (B)~Banks and Bihari's method,
    (G)~Algorithm~\ref{alg:grd_est}.
  }
  \label{fig:exp2_stats}
  \end{minipage}%
  \hspace{2em}%
  \begin{minipage}[c]{.3\linewidth}
  \includegraphics[width=\linewidth,trim=5 5 14 10,clip]{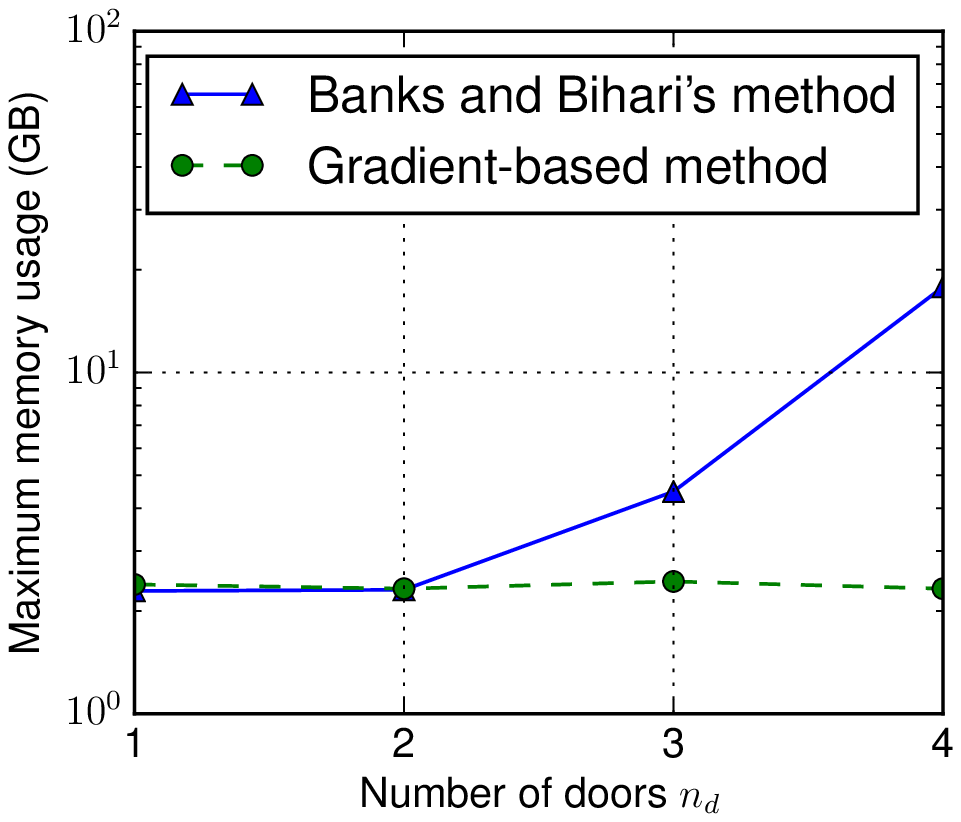}
  \caption{
    Memory usage in $\unit{GB}$ of the experiments in Sec.~\ref{subsec:exp3}.
  }
  \label{fig:mem}
  \end{minipage}
\end{figure}

Figures~\ref{fig:exp2_dr} and~\ref{fig:exp2_tem} are analogous to those in Figure~\ref{fig:exp1_stats}, while Figure~\ref{fig:exp2_dr1} shows the estimation error just for door $d_1$, which is located very close to the thermostat $s_1$.
As shown in these figures, both estimation algorithms do an almost equally poor job at estimating the door configuration, and Algorithm~\ref{alg:grd_est} is marginally better at estimating the initial temperature distribution.
Yet, both algorithms are capable of accurately estimating the configuration of the door closest to the thermostat, which shows that even with one thermostat we can get enough information of the resident's activities.

\subsection{Memory usage comparison}
\label{subsec:exp3}

A significant advantage of Algorithm~\ref{alg:grd_est} when compared to probabilistic estimation algorithms is that our method does not need to compute numerical solutions of the set of differential equations for each possible configuration $\theta \in \set{0,1}^{n_d}$.
Figure~\ref{fig:mem} shows the maximum memory usage of both algorithm implementations as the number of doors to estimate increases from~1 to~4.
Our results show that the probabilistic estimation method can be used only for small values of $n_d$, quickly outgrowing the amount of memory in standard computers (for $n_d=4$ the usage was $16 \unit{GB}$ approx.), while Algorithm~\ref{alg:grd_est} memory usage remains almost constant (at $2.4 \unit{GB}$ approx.).

%% file: sections/conclusion.tex
\section{CONCLUSION}
\label{sec:conclusion}

Our gradient-based estimation method and simulation results show the potential for reconstructing indoor climate and building configuration by using only thermostat sensor data, thus reducing the need for extra sensors to monitor a smart buildings.
Also, since the method can accurately estimate the indoor climate and configuration with acceptable memory usage, it can be used in coordination to advanced MPC control strategies, significantly increasing the efficiency of HVAC units without a decrease in human comfort.
Our method has the potential to enable interesting new applications.
For example, since it is able to identify a building's configuration in real-time, it can potentially be applied to monitor an unexpected break-in.

%% file: sections/appendx.tex
\section{DERIVATION OF ADJOINT EQUATIONS}
\label{sec:appendix_a}

Consider the Lagrangian function in equation~\eqref{eq:lg}.
For each set of functions $(w, v, q)$ in the respective dual space of the tuple $( Te, u, p )$, we can write $\dprodb{\pderiv{L}{T_e}}{w}_{\OM \times [0,T]} = 0$ as follows:
\begin{equation}
  \label{eq:dir_w1}
    \dprodB{\pderiv{L}{T_e}}{w}_{\OM \times [0,T]}
    = \dprodB{\pderiv{J}{T_e}}{w}_{ \OM \times \sqparen{0,T} }
    + \dprodB{\lambda_1}{
      \pderiv{w}{t}
      - \grad[x] \cdot \p{\kappa\, \grad[x]{w}}
      + u \cdot \grad[x]{w} }_{\OM \times [0,T]}
    + \dprod{\lambda_4}{w}_{\partial \OM \times [0,T]}
    + \dprod{\lambda_6}{w(0,\cdot) }_{\OM}
    = 0.
\end{equation}
Similarly, we can write $\dprodb{\pderiv{L}{u}}{v}_{\OM} = 0$ as:
\begin{equation}
    \dprodB{\pderiv{L}{u}}{v}_{\OM}
    =
    \dprodB{\lambda_1}{
      v \cdot \grad[x]{T_e} }_{\OM \times [0,T]}
    + \dprodB{\lambda_2}{
      - \frac{1}{\reno} \lapl[x]{v}
      + \p{u \cdot \grad[x]{}} v
      + \p{v \cdot \grad[x]{}} u
      + \alpha\, v }_{\OM}
    + \dprod{\lambda_3}{\diver[x]{v}}_{\OM}
    + \dprod{\lambda_5}{v}_{\Gamma_w}
    = 0,
  \label{eq:dir_v}
\end{equation}
and we can write $\dprodb{\pderiv{L}{p}}{q}_{\OM} = 0$ as:
$\dprodb{\lambda_2}{\grad[x]{q}}_{\OM} = 0$.

Applying integration by parts and Green's formula, equations~\eqref{eq:dir_w1} and~\eqref{eq:dir_v} become:
\begin{multline}
  \label{eq:dir_w2}
    \dprodB{\pderiv{J}{T_e}}{w}_{\OM \times \sqparen{0,T}}
    - \dprodB{
      \pderiv{\lambda_1}{t}
      + \grad[x] \cdot \p{ \kappa\, \grad[x]{ \lambda_1 } }
      + u \cdot \grad[x]{\lambda_1} } {w}_{ \OM \times \sqparen{0,T} }
    + \dprodb{ \lambda_1(\cdot,T) }{ w\paren{\cdot,T} }_{\OM} + \\
    + \dprodb{ \lambda_6-\lambda_1(\cdot,0)}{ w(\cdot,0) }_{\OM}
    + \dprodB{
      \kappa\, \pderiv{\lambda_1}{\vec{n}}
      + \lambda_4
      + \vec{n} \cdot u \lambda_1}{w}_{\partial \OM \times \sqparen{0,T} }
    - \dprodB{\kappa\, \lambda_1}{\pderiv{w}{\vec{n}}}_{\partial \OM \times \sqparen{0,T}}
    = 0,
\end{multline}
\begin{multline}
  \label{eq:dir_v2}
  \dprodb{\lambda_1}{
    v \cdot \grad[x]{T_e} }_{\OM \times [0,T]}
  + \dprodB{
    \alpha\, \lambda_2
    - \frac{1}{\reno} \lapl[x]{\lambda_2}
    + \grad[x]{u} \cdot \lambda_2
    - u \cdot \grad[x]{\lambda_2}
    - \grad[x]{\lambda_3} } {v}_{\OM} + \\
  - \frac{1}{\reno}
  \dprodb{\lambda_2}{\pderiv{v}{\vec{n}}}_{\partial \OM}
  + \dprodb{
    \frac{1}{\reno} \pderiv{\lambda_2}{\vec{n}}
    + \lambda_3\, \vec{n}
    + \paren{ u \cdot \vec{n} } \lambda_2 }{ v }_{\partial \OM}
  + \dprodb{ \lambda_5 }{ v }_{\Gamma_w}
  = 0,
\end{multline}
where $\vec{n}$ is the vector normal to the boundary at $x \in \bdry{\Omega}$.

From the identities above it follows that, in order to make $\dprodb{\pderiv{L}{T_e}}{w}_{\OM \times [0,T]} = 0$, $\dprodb{\pderiv{L}{u}}{v}_{\OM} = 0$, and $\dprodb{\pderiv{L}{p}}{q}_{\OM} = 0$ for any set of functions $(w,v,q)$, a sufficient condition for the dual variables $\lambda_{1,2,3,6}$ is to satisfy the differential equations~\eqref{eq:ad_pde1} to~\eqref{eq:ad_pde3} and their boundary conditions.

\section{DERIVATION OF FR\'ECHET DERIVATIVES}
\label{sec:appendix_b}

As explained in Section~\ref{sec:adj_eq}, if we take variations $(\delta\pi_0, \delta\theta)$ of our optimization variables, those will induce variations $\delta\alpha$, $\delta\kappa$, $\delta T_e$, $\delta u$, and $\delta p$.
Then, from equations~\eqref{eq:temp},~\eqref{eq:ns1} and~\eqref{eq:ns2}, it follows that the variations satisfy the following differential equations:
\begin{gather}
  \label{eq:pde_dlt1}
  \pderiv{\delta T_e}{t}
  - \grad[x]{} \cdot \paren{ \delta \kappa\, \grad[x]{T_e} }
  - \grad[x]{} \cdot \paren{ \kappa\, \grad[x]{\delta T_e} }
  + \delta u \cdot \grad[x]{T_e}
  + u \cdot \grad[x]{\delta T_e} = 0,
  \\
  \label{eq:pde_dlt2}
  \delta \alpha\, u
  + \alpha\, \delta u
  - \frac{1}{\reno}\, \lapl[x]{\delta u}
  + \delta u \cdot \grad[x]{u}
  + u \cdot \grad[x]{\delta u}
  + \grad{\delta p} = 0,
  \\
  \label{eq:pde_dlt3}
  \diver[x]{\delta p} = 0,
\end{gather}
with the following boundary and initial conditions: $\delta T_e(x,t) = 0$ for each $x \in \bdry{\Omega}$ and $t \in [0,T]$, $\delta T_e(x,0) = \delta\pi_0(x)$ for each $x \in \Omega$, and $\delta u(x) = 0$ for each $x \in \Gamma_w$.

Now, using equations~\eqref{eq:cost}, \eqref{eq:pde_dlt1}, \eqref{eq:pde_dlt2} and~\eqref{eq:pde_dlt3} and their boundary and initial conditions, we get:
\begin{multline}
  \label{eq:cost_dif_0}
  J\paren{ \theta + \delta \theta, \pi_0 + \delta \pi_0}
  - J\paren{ \theta, \pi_0 }
  =
  \dprodB{\pderiv{J}{T_e}}{\delta T_e}_{\OM \times \sqparen{0,T}}
  + \dprodB{\pderiv{J}{T_e}}{\delta \pi_0}_{\OM} + \\
  + \dprodB{\lambda_1}{
    \pderiv{\delta T_e}{t}
    - \grad[x]{} \cdot \paren{ \delta \kappa\, \grad[x]{T_e} }
    - \grad[x]{} \cdot \paren{ \kappa\, \grad[x]{\delta T_e} }
    + \delta u \cdot \grad[x]{T_e}
    + u \cdot \grad[x]{\delta T_e}}_{ \OM \times \sqparen{0,T} } + \\
  + \dprodB{\lambda_2} {
    \delta \alpha\, u
    + \alpha\, \delta u
    - \frac{1}{\reno} \lapl[x]{\delta u}
    + \delta u \cdot \grad[x]{u}
    + u \cdot \grad[x]{\delta u}
    + \grad[x]{\delta p} }_{\OM}
  + \dprodb{\lambda_3}{\diver[x]{\delta u}}_{\OM} + \\
  + \dprodb{\lambda_4}{\delta T_e}_{ \partial \OM \times \sqparen{0,T} }
  + \dprodb{\lambda_5}{\delta u}_{\Gamma_w}
  + \dprodb{\lambda_6}{\delta Te(\cdot,0) - \delta \pi_0}_{\OM}
\end{multline}
where $\set{\lambda_i}_{i=1}^6$ are the adjoint variables defined in Section~\ref{sec:adj_eq}.
Then, applying integration by parts and Green's formula equation~\eqref{eq:cost_dif_0}, and after canceling terms using the identities in equations~\eqref{eq:ad_pde1} to~\eqref{eq:ad_pde3}, we can get:
\begin{equation}
  \label{eq:cost_diff2}
    J\paren{ \theta + \delta \theta, \pi_0 + \delta \pi_0 }
    - J\paren{ \theta, \pi_0 }
    =
    \dprodb{ \grad[\pi_0]{J} - \lambda_6 }{\delta \pi_0}_{\OM}
    + \dprodb{ \grad[x]{ \lambda_1 } \cdot \grad[x]{T_e} }{ \delta \kappa }_{\OM \times [0,T]}
    + \dprodb{ \lambda_2 \cdot u }{ \delta \alpha }_{\OM},
\end{equation}
which are equivalent to the directional derivatives in equations~\eqref{eq:alpha_kappa_deriv} and~\eqref{eq:theta_i0_deriv}, as desired.